\documentclass[12pt]{amsart}
\usepackage[latin1]{inputenc}
\usepackage{tikz}
\allowdisplaybreaks

\newtheorem{theorem}{Theorem}

\title[On protected nodes in Digital Search Trees]
{On protected nodes in Digital Search Trees}

\author[R.R.X.~DU]{Rosena R.X. Du}

\address{Rosena R.X. Du, Department of Mathematics, East China Normal University\\
500 Dongchuan Road, Shanghai, 200241, P. R. China.}
\email{rxdu@math.ecnu.edu.cn}

\author[H.~Prodinger]{Helmut Prodinger}

\address{Helmut Prodinger,
Mathematics Department, Stellenbosch University,
7602 Stellenbosch, South Africa.}
\email{hproding@sun.ac.za}


\begin{document}

\dedicatory{Dedicated to Philippe Flajolet (1948--2011)}

\begin{abstract}Recently, 2-protected nodes were studied in the context of ordered trees and $k$-trees.
These nodes have a distance of at least 2 to each leaf. Here, we study digital search trees, which are
binary trees, but with a different probability distribution underlying. Our result says, that \emph{grosso modo}
some $31\%$ of the nodes are 2-protected. Methods include exponential generating functions, contour integration,
and some elements from $q$-analysis.
\end{abstract}

\maketitle

\section{Introduction}

Cheon and Shapiro~\cite{ChSh08} started the study of 2-protected nodes in trees. A node enjoys this
property if its distance to any leaf is at least 2. A simpler notion is 1-protected: exactly the nodes that
are not leaves are 1-protected. In the cited paper, the family of ordered trees was considered, and it
was found that asymptotically a proportion of $\frac16$ of the nodes is 2-protected. Recently, Mansour~\cite{Mansour11}
complemented these results by studying $k$-ary trees.

In the present note, we study the analogous quantity for \emph{Digital Search Trees} (DSTs), a structure that
is important in Computer Science~\cite{Knuth98}. As trees, they are binary trees, but the (probability) distribution
is quite different. From a mathematical point of view, they always lead to interesting and nontrivial considerations,
with a flair of $q$-analysis. Here are a few papers of relevance:~\cite{FlSe86, Prodinger92, KiPr88, LoPr06, HwFuZa10}

\begin{figure}[h!]
\begin{center}
\begin{tikzpicture}[
    s1/.style={
    circle,
    very thick,
    draw=black,
    top color=white,
    bottom color=yellow!50!black!30
},
    s2/.style={
    circle,
    draw=white!50!black!80,
    top color=white,
    bottom color=yellow!50!black!30
}]
\path node at (1,0) {$A:1001$};%
\path node at (1,-.5) {$B:0110$};%
\path node at (1,-1) {$C:0000$};%
\path node at (1,-1.5) {$D:1111$};%
\path node at (1,-2) {$E:0100$};%
\path node at (1,-2.5) {$F:0101$};%
\path node at (1,-3) {$G:1101$};%
\path node at (1,-3.5) {$H:1110$};%
\path node at (1,-4) {$I:1100$};%
\footnotesize %
\path node(a) at (7,0) [s1, text centered]{$A$};%
\path node(b) at (5,-1.5) [s2, text centered]{$B$};%
\path node(c) at (4,-3) [s2, text centered]{$C$};%
\path node(d) at (9,-1.5) [s1, text centered]{$D$};%
\path node(e) at (6,-3) [s2, text centered]{$E$};%
\path node(f) at (5,-4.5) [s2, text centered]{$F$};%
\path node(g) at (10,-3) [s2, text centered]{$G$};%
\path node(h) at (11,-4.5) [s2, text centered]{$H$};%
\path node(i) at (9,-4.5) [s2, text centered]{$I$};%
\path [draw,-,black!90] (a) -- (b) node[above,pos=.6,black]{{\tiny $0$}};%
\path [draw,-,black!90] (b) -- (c) node[above,pos=.6,black]{{\tiny $0$}};%
\path [draw,-,black!90] (a) -- (d) node[above,pos=.6,black]{{\tiny $1$}};%
\path [draw,-,black!90] (b) -- (e) node[above,pos=.6,black]{{\tiny $1$}};%
\path [draw,-,black!90] (e) -- (f) node[above,pos=.6,black]{{\tiny $0$}};%
\path [draw,-,black!90] (d) -- (g) node[above,pos=.6,black]{{\tiny $1$}};%
\path [draw,-,black!90] (g) -- (h) node[above,pos=.6,black]{{\tiny $1$}};
\path [draw,-,black!90] (g) -- (i) node[above,pos=.6,black]{{\tiny $0$}};%
\end{tikzpicture}
\end{center}
\caption{\emph{A digital search tree with nine nodes, among which $A$ and $D$ are 2-protected.}}
\label{fg-dst}
\end{figure}
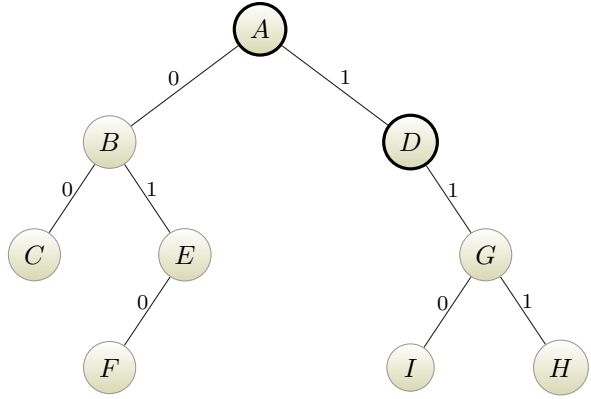

DSTs are constructed as follows. Given a sequence of binary strings, we place the first in the root node; those starting with ``0" (``1") are directed to the left (right) subtree of the root, and are constructed recursively by the same procedure but with the removal of their first bits when comparisons are made. See Figure \ref{fg-dst} for an illustration.

In the following section we will show that the proportion of 2-protected nodes in the DST model is about $31\%$; a more detailed
statement will be given later.

We collect here are few notations. These quantities belong to the realm of $q$-series and can be found in~\cite{Andrews76}, although with a slightly different notation:
\begin{equation*}
Q_m=\prod_{k=1}^m\Big(1-\frac1{2^k}\Big),
\quad Q_\infty=\prod_{k=1}^\infty\Big(1-\frac1{2^k}\Big),
\quad Q(x)=\prod_{k=1}^\infty\Big(1-\frac x{2^k}\Big).
\end{equation*}
There is a formula that is equivalent to one of Euler's partition identities:
\begin{equation*}
Q(t)=\sum_{m\ge0}a_{m+1}t^m \qquad\text{with}\qquad
\quad a_{m+1}=\frac{(-1)^m2^{-\binom{m+1}2}}{Q_m}.
\end{equation*}
Finally, we will use $L=\log2$.

\section{Average number of 2-protected nodes}

Denote by $l_n$ the average number of 2-protected nodes in a random DST, built from $n$ data. By random we mean that whenever a decision
has to be made whether to go down to the left or right, a fair coin is tossed, and a direction is chosen with probability $\frac12$.

The following recursion follows from the observation that, provided we have $n+1$ data, $k$ go to the left and $n-k$ go to the
right, and such a split happens with probability $\binom nk2^{-n}$. One node goes to the root and is always 2-protected except in
the instances $k=1$ or $k=n-1$. Therefore
\begin{align*}
l_{n+1}&=\sum_{k=0}^n\binom nk2^{-n}\Big(l_k+l_{n-k}+1\Big)
-\sum_{k= 1\text{ or }  n-1}\binom nk2^{-n}\\
&=1+2^{1-n}\sum_{k=0}^n\binom nkl_k-n2^{1-n}.
\end{align*}
This recursion is true for $n\ge3$, with initial conditions $l_0=l_1=l_2=0$, $l_3=\frac12$.
Our treatment follows \cite{FlSe86}. We introduce the exponential generating function
$L(z)=\sum_{n\ge0}l_nz^n/n!$ and translate the recursion:
\begin{equation*}
\sum_{n\ge3}l_{n+1}\frac{z^n}{n!}=\sum_{n\ge3}\frac{z^n}{n!}+\sum_{n\ge3}\frac{z^n}{n!}\sum_{k=0}^n\binom nk2^{1-n}l_k
-\sum_{n\ge3}n2^{1-n}\frac{z^n}{n!}
\end{equation*}
or
\begin{equation*}
\sum_{n\ge0}l_{n+1}\frac{z^n}{n!}-l_{3}\frac{z^2}{2!}=\sum_{n\ge3}\frac{z^n}{n!}+\sum_{n\ge0}\frac{z^n}{n!}\sum_{k=0}^n\binom nk2^{1-n}l_k
-\sum_{n\ge3}n2^{1-n}\frac{z^n}{n!},
\end{equation*}
which leads after some simple manipulations to
\begin{equation*}
L'(z)=e^z-ze^{z/2}-1+\frac{z^2}{4}+2e^{z/2}L(\frac z2).
\end{equation*}
Now we introduce the \emph{Poisson generating function} $M(z)=e^{-z}L(z)=\sum_{n\ge0}m_nz^n/n!$ and rewrite the equation:
\begin{equation*}
M'(z)+M(z)=1-ze^{-z/2}-e^{-z}+\frac{z^2}{4}e^{-z}+2M(\frac z2).
\end{equation*}
For $n\ge1$, we can read off the coefficients of $z^n/n!$:
\begin{equation*}
m_{n+1}=-(1-2^{1-n})m_n+n(-1)^{n}2^{1-n}-(-1)^n+\frac{n(n-1)}{4}(-1)^n.
\end{equation*}
In order to solve it, we rewrite it as
\begin{equation*}
\frac{m_{n+1}(-1)^n}{Q_{n-1}}=\frac{m_n(-1)^{n-1}}{Q_{n-2}}+\frac{n2^{1-n}-1+\frac{n(n-1)}{4}}{Q_{n-1}},
\end{equation*}
which can be summed and leads to
\begin{equation*}
\frac{m_{N+1}(-1)^N}{Q_{N-1}}=\sum_{n=2}^N\frac{n2^{1-n}-1+\frac{n(n-1)}{4}}{Q_{n-1}}
\end{equation*}
and eventually to
\begin{equation*}
m_{N}=Q_{N-2}(-1)^N\sum_{n=1}^{N-2}\frac{1-(n+1)2^{-n}-\frac{n(n+1)}{4}}{Q_{n}}.
\end{equation*}
Since
\begin{equation*}
l_n=\sum_{k=2}^n\binom nk m_k
\end{equation*}
we found the following explicit formula that we formulate as a theorem.
\begin{theorem}
The average number of 2-protected nodes in random DSTs of size $N\ge1$ is exactly given by
\begin{equation*}
l_ N=\sum_{k=2}^N\binom Nk (-1)^k
Q_{k-2}\sum_{n=1}^{k-2}\frac{1-(n+1)2^{-n}-\frac{n(n+1)}{4}}{Q_{n}}.
\end{equation*}
\end{theorem}

Now we turn to the asymptotic evaluation of $l_N$ as $N\to\infty$. Again, we follow the approach in \cite{FlSe86}
and use Rice's integrals, which means that we are able to rewrite $l_N$ as a contour integral. Changing the contour
of integration and collecting residues produces the asymptotic expansion of interest. Many examples have been described
in \cite{FlSe95}. In order to do so, one must extend the function
\begin{equation*}
Q_{k-2}\sum_{n=1}^{k-2}\frac{1-(n+1)2^{-n}-\frac{n(n+1)}{4}}{Q_{n}}
\end{equation*}
so that it makes sense for any complex $k$, not just integers. This will be discussed now.

We have $Q_{k-2}=Q_\infty/Q(2^{2-k})$, and this makes sense for any $k$.
Now we have, using Euler's identity mentioned in the Introduction,
\begin{equation*}
\frac1{Q_n}=\frac{Q(2^{-n})}{Q_\infty}=\frac{1}{Q_\infty}\sum_{m\ge0}a_{m+1}2^{-nm},
\end{equation*}
and this makes sense for any $n$, since the smallness of the $a_m$'s handles all convergence issues.
Therefore
\begin{equation*}
\sum_{n=1}^{k-2}\frac{1-(n+1)2^{-n}-\frac{n(n+1)}{4}}{Q_{n}}=
\frac{1}{Q_\infty}\sum_{m\ge0}a_{m+1}\sum_{n=1}^{k-2}\Big[1-(n+1)2^{-n}-\frac{n(n+1)}{4}\Big]2^{-nm}.
\end{equation*}
The inner sum (on $n$) can be explicitly evaluated, but since it is long and ugly, we don't display it
here. The resulting form (that we keep in our Maple calculation) can be used for any $k\in\mathbb{C}$.

The integral expression is
\begin{equation*}
l_N=-\frac1{2\pi i}\int_{\mathcal{C}}\frac{\Gamma(N+1)\Gamma(-z)}{\Gamma(N+1-z)}\psi(z)dz,
\end{equation*}
where $\mathcal{C}$ encircles the poles $2,3,\dots,N$ and no others. The function $\psi(z)$ is the extension of
\begin{equation*}
Q_{k-2}\sum_{n=1}^{k-2}\frac{1-(n+1)2^{-n}-\frac{n(n+1)}{4}}{Q_{n}}
\end{equation*}
as just discussed. Changing the contour, one encounters other poles. They must be subtracted and produce the
asymptotic expansion that we need. The main contribution comes from $z=1$. There are also poles at $z=1+\chi_k$,
with $\chi_k=\frac{2\pi i k}{L}$, and they contribute a tiny oscillating function $N\cdot \delta(\log_2N)$, where
the amplitude of $\delta(x)$ is typically smaller than $10^{-5}$. In order to keep this note short and crisp, we refrain
from computing this function explicitly. It is not difficult, and there are many similar examples in the literature.
So we concentrate now on $z=1$, and we will find a simple pole. As a first step, we consider
\begin{gather*}
\lim_{k\to1}\frac{Q_{k-1}}{1-2^{1-k}}\sum_{n=1}^{k-2}\Big[1-(n+1)2^{-n}-\frac{n(n+1)}{4}\Big]2^{-nm}.
\end{gather*}
This limit can be computed by Maple, with the result
\begin{equation*}
b_m:=\frac{1}{4L}{\frac {B(2^{-m})}{ ( 2^{-m}-1 )^3(2^{-m}-2)^2 }}
\end{equation*}
and $B(x):=16L-48xL+48{x}^{2}L-16{x}^{3}L -20
+60x -69{x}^{2} +36{x}^{3}-7{x}^{4}-8x\log  \left( x \right) +
12{x}^{2}\log  \left( x \right) -10{x}^{3}\log  \left( x \right) +4
{x}^{4}\log  \left( x \right) $. Note that $b_0$ is interpreted as a limit:
\begin{equation*}
b_0= \frac{37}{12L}-4.
\end{equation*}
So we are left with the negative residue of
\begin{equation*}
-\frac{\Gamma(N+1)\Gamma(-z)}{\Gamma(N+1-z)}
\end{equation*}
at $z=1$, which is just $N$. Summarizing, we found the asymptotic behaviour.
\begin{theorem}
The average number $l_N$ of 2-protected nodes in random DSTs of size $N$ admits the asymptotic expansion
\begin{equation*}
l_N= N\cdot \frac{1}{Q_\infty}\sum_{m\ge0}a_{m+1}b_m+N\cdot\delta(\log_2N)+O(1),
\end{equation*}
where the numerical constant evaluates to
$0.30707981393605921828549\dots$. The tiny periodic function $\delta(x)$ has a Fourier expansion that could be computed
in principle. The remainder term $O(1)$ stems from the next pole at $z=0$.
\end{theorem}

One referee has suggested to give the explicit expression of the periodic function $\delta(x)$ without proof. Here it is:
\begin{align*}
\delta(x)&=\frac1{Q_{\infty}}\sum_{l\neq0}\sum_{m\ge0}a_{m+1} \frac{l\pi 2^m}
{2L^2(2^m-1)^2(2^{m+1}-1)}\\
&\qquad\qquad\times\Big[iL(7-15\cdot2^m+10\cdot 4^m)-2\pi l(2^{m+1}-1)\Big]e^{-2\pi i lx}.
\end{align*}

For example, $l_{500}/500=0.305710\dots$.

\subsection*{Remark}Flajolet and Sedgewick in \cite{FlSe86} solved an open problem of Knuth~\cite{Knuth98}, and considered the number of endnodes.
They found this to be on average as $\beta\cdot N$, with $\beta=0.372046812\dots$. Again, there are tiny oscillations.
The quantity $(1-\beta)N$ is (asymptotically) the number of 1-protected nodes. So, there are roughly $63\%$ 1-protected nodes,
and our new results say that there are about $31\%$ 2-protected nodes.

\subsection*{Acknowledgement} The first author is partially supported by the National Science Foundation of China under Grant 10801053, and the Shanghai Rising-Star Program (No. 10QA1401900). The second author is supported by an International Science and Technology Agreement (Grant 67215) from the NRF (South Africa).


\begin{thebibliography}{10}

\bibitem{Andrews76}
George~E. Andrews.
\newblock {\em The theory of partitions}.
\newblock Addison-Wesley Publishing Co., Reading, Mass.-London-Amsterdam, 1976.
\newblock Encyclopedia of Mathematics and its Applications, Vol. 2.

\bibitem{ChSh08}
Gi-Sang Cheon and Louis~W. Shapiro.
\newblock Protected points in ordered trees.
\newblock {\em Appl. Math. Lett.}, 21(5):516--520, 2008.

\bibitem{FlSe86}
Philippe Flajolet and Robert Sedgewick.
\newblock Digital search trees revisited.
\newblock {\em SIAM J. Comput.}, 15(3):748--767, 1986.

\bibitem{FlSe95}
Philippe Flajolet and Robert Sedgewick.
\newblock Mellin transforms and asymptotics: finite differences and {R}ice's
  integrals.
\newblock {\em Theoret. Comput. Sci.}, 144(1-2):101--124, 1995.
\newblock Special volume on mathematical analysis of algorithms.

\bibitem{HwFuZa10}
Hsien-Kuei Hwang, Michael Fuchs, and Vytas Zacharovas.
\newblock Asymptotic variance of random symmetric digital search trees.
\newblock {\em Discrete Math. Theor. Comput. Sci.}, 12(2):103--165, 2010.

\bibitem{KiPr88}
Peter Kirschenhofer and Helmut Prodinger.
\newblock Eine {A}nwendung der {T}heorie der {M}odulfunktionen in der
  {I}nformatik.
\newblock {\em \"Osterreich. Akad. Wiss. Math.-Natur. Kl. Sitzungsber. II},
  197(4-7):339--366, 1988.

\bibitem{Knuth98}
Donald~E. Knuth.
\newblock {\em The Art of Computer Programming}, volume 3: Sorting and
  Searching.
\newblock Addison-Wesley, 1973.
\newblock Second edition, 1998.

\bibitem{LoPr06}
Guy Louchard and Helmut Prodinger.
\newblock Asymptotics of the moments of extreme-value related distribution
  functions.
\newblock {\em Algorithmica}, 46(3-4):431--467, 2006.

\bibitem{Mansour11}
Toufik Mansour.
\newblock Protected points in {$k$}-ary trees.
\newblock {\em Appl. Math. Lett.}, 24(4):478--480, 2011.

\bibitem{Prodinger92}
Helmut Prodinger.
\newblock External internal nodes in digital search trees via {M}ellin
  transforms.
\newblock {\em SIAM J. Comput.}, 21(6):1180--1183, 1992.

\end{thebibliography}
\end{document}